\documentclass[12 pt]{amsart}
\usepackage{amsmath}
\usepackage{amssymb}
\usepackage{color}
\usepackage[dvips]{epsfig}

\def\<{\left < }
\def\>{\right >}
\def\({\left ( }
\def\){\right )}

\def\a{\alpha}
\def\b{\beta}
\def\g{\gamma}
\def\t{\theta}
\def\l{\lambda}
\def\m{\mu}
\begin{document}

\title [Cheng-Yau operator and   Gauss map of surfaces of revolution]
{Cheng-Yau operator and   Gauss map of surfaces of revolution}

\author{Dong-Soo Kim$^1$, Jong Ryul Kim and Young Ho Kim$^3$}

\address{\newline Department of Mathematics, Chonnam National University, Gwangju
500-757,  South Korea \\ \newline
Department of Mathematics, Kunsan National University,
Kunsan 573-701, South Korea
\\ \newline
Department of Mathematics, Kyungpook
National University, Daegu 702-701,  South Korea}
 \email{
dosokim@chonnam.ac.kr, kimjr0@yahoo.com and yhkim@knu.ac.kr}

\thanks{\newline\indent { $^1$  supported by Basic Science Research Program through the National Research Foundation of Korea (NRF) funded by the Ministry of Education, Science and Technology (2010-0022926).}
\newline\indent { $^3$   supported by Basic Science Research Program through the National Research Foundation of Korea (NRF)
funded by the Ministry of Education, Science and Technology (2012R1A1A2042298) and supported by Kyungpook
National University Research Fund, 2012.}}
\begin{abstract}
We study  the Gauss map $G$ of surfaces of revolution
in the 3-dimensional Euclidean space ${\mathbb{E}^3}$
with respect to the so called Cheng-Yau operator $\square$
 acting on the functions defined on the surfaces.
As a result, we establish the classification theorem  that the only
 surfaces of revolution with  Gauss map  $G$ satisfying $\square G=AG$
 for some $3\times3$ matrix $A$ are the
planes, right circular cones, circular cylinders and spheres.
\end{abstract}
\vskip 1cm

\maketitle

\section{Introduction }

The theory of Gauss map of a surface in a Euclidean space and a pseudo-Euclidean space
 is always one of interesting topics
 and it has been investigated
from the various viewpoints by many differential geometers
(\cite{bb, bv, cho1, cho2, ckky, cp, dpv, du, kkkr, k1, k2, kky, ks, kt, ky}).

  Let  $M $ be a surface of the Euclidean 3-space
$\mathbb{E}^3$.
The map
$G:M \rightarrow S^{2} \subset {\mathbb{E}^3}$ which sends each
point of $M$ to the unit normal vector to $M$ at the point is
called the
  {\itshape Gauss map} of the
surface
$M,$ where
$S^{2}$ is the unit sphere in $\mathbb{E}^3$ centered at
the origin. It is well known that $M$ has constant mean curvature if and
only if $\Delta G = ||dG||^2G$, where $\Delta$ is the Laplace operator on M corresponding
to the induced metric on M from $\mathbb{E}^3$ (\cite{rv}).
Surfaces whose Gauss map is an eigenfunction of  Laplacian, that is, $\Delta G = \l G$
for some constant $\l\in R$, are the planes,  circular cylinders and spheres (\cite{cp}).

Generalizing  this equation, F. Dillen, J. Pas and
L. Verstraelen (\cite{dpv}) studied surfaces of revolution in a Euclidean 3-space $\mathbb{E}^3$
such that its Gauss map G satisfies the condition
\begin{equation}
  \tag{1.1}
  \Delta G=AG, \quad A\in R^{3\times3}.
\end{equation}

 As a result, they proved (\cite{dpv})

\noindent {\bf Proposition 1.1.} Among the surfaces of revolution in ${\mathbb{E}^3}$,
 the only ones whose Gauss map satisfies (1.1) are the planes, the spheres and the circular cylinders.

Baikoussis and Blair also studied ruled surfaces in ${\mathbb{E}^3}$ and proved (\cite{bb})

\noindent {\bf Proposition 1.2.} Among the ruled surfaces in ${\mathbb{E}^3}$,
 the only ones whose Gauss map satisfies
(1.1) are the planes and the circular cylinders.

Generalized slant cylindrical surfaces (GSCS's) are natural extended notion of
surfaces of revolution (\cite{kk}).
Surfaces of revolution, cylindrical surfaces and tubes along a plane curve
  are special cases of  GSCS's. In \cite{ks}, the first author and B. Song
proved that among the GSCS's in ${\mathbb{E}^3}$,
 the only ones whose Gauss map satisfies (1.1) are the planes, the spheres and the circular cylinders.

The so-called Cheng-Yau operator $\square$ (or, $L_1$) is a natural extension of the Laplace operator
 $\Delta$ (cf. \cite{ag}, \cite{cy}). Hence, following the condition (1.1),
  it is natural to ask as follows.

 \noindent {\bf Question 1.3.}
Among the  surfaces of revolution  in a
Euclidean 3-space $\mathbb{E}^3$, which one satisfy the  following condition?
\begin{equation}
  \tag{1.2}
  \square G=AG,\quad A\in R^{3\times3}.
\end{equation}

In this paper, we give a complete answer to the above question.

Throughout this paper, we assume that all objects are smooth and connected,
 unless otherwise mentioned.
\vskip0.5cm

\section{Cheng-Yau operator and Lemmas }
\vskip0.5cm

Let $M$ be an oriented surface in $E^3$ with Gauss map $G$.
We denote by $S$ the shape operator of $M$ with respect to the Gauss map $G$.
For each $k = 0, 1$, we put $P_0 = I, P_1 = tr(S)I - S$, where
$I$ is the identity operator acting on the tangent bundle of M.
Let us define an operator
$L_k : C^{\infty}(M) \rightarrow  C^{\infty}(M)$ by
 $L_k(f) = -tr(P_k \circ \nabla ^2f)$,
 where $\nabla ^2f: \chi(M) \rightarrow  \chi(M)$ denotes the self-adjoint linear operator
 metrically equivalent to the hessian of $f$.
 Then, up to signature,  $L_k$ is the
linearized operator of the first variation of the $(k+1)$-th mean curvature arising
from normal variations of the surface.
Note that the
operator $L_0$ is nothing but the Laplace operator acting on $M$, i.e., $L_0 = \Delta$
and $L_1 = \square$ is called the Cheng-Yau operator introduced in \cite{cy}.
\vskip0.3cm

Now, we state a useful lemma as follows (\cite{ag}).
\vskip0.3cm

\noindent {\bf Lemma 2.1.} Let $M$ be an oriented surface in $E^3$ with Gaussian curvature $K$
and mean curvature $H$. Then, the Gauss map $G$ of $M$ satisfies
\begin{equation}
\tag{2.1} \square G = \nabla K + 2HK G,
\end{equation}
where $\nabla K$ denotes the gradient of $K$.
\vskip0.3cm

Now, using Lemma 2.1  we give  some examples of surfaces  with
Gauss map satisfying (1.2).
\par \vskip0.3cm

 \noindent {\bf Examples.}

 \noindent (1) Flat surfaces.  In this case, we have $\square G=0$,  and hence
 flat surfaces  satisfy $\square G=AG$  for some  $3\times 3 $ matrix $A$.
 Note that the matrix $A$ must be singular.

  \noindent (2) Spheres: $(x-a)^2+(y-b)^2+(z-c)^2=r^2$.  In this case, we have
 $G=\frac{1}{r}(x-a,y-b,z-c)$ so
 the sphere  satisfies $\square G=AG$  with $A=\frac{-1}{r^3}I$, where $I$ denotes the identity matrix.

\vskip.8cm

\section{Gauss map of surfaces of revolution}

\vskip.8cm

We consider  a unit speed plane curve  $C:(x(s), 0, z(s))$ with $x(s)>0$ in the $xz$ plane
which is defined on an interval $I$. By rotating the curve $C$ around $z$-axis,
we get a surface of revolution $M$, which is parametrized by
\begin{equation}
\tag{3.1}
X(s,t)=(x(s)\cos t, x(s)\sin t, z(s)).
\end{equation}

The adapted  frame field $\{e_1, e_2, G\}$ on the surface of revolution $M$ are given by
\begin{equation}
\tag{3.2}
\begin{aligned}
e_1&=X_s=(x'(s)\cos t, x'(s)\sin t, z'(s)),\\
 e_2&=\frac{1}{x}X_t=(-\sin t, \cos t, 0), \\
 G&=e_1\times e_2=(-z'\cos t,-z'\sin t, x').
\end{aligned}
\end{equation}
The principal curvatures $k_1, k_2$ of $M$ with respect to the Gauss map $G$ are respectively
(\cite{do})
\begin{equation}
\tag{3.3}
\begin{aligned}
k_1&=\<S(e_1),e_1\>=x'z''-x''z'=\kappa,\\
k_2&=\<S(e_2),e_2\>=\frac{z'}{x},
\end{aligned}
\end{equation}
where $S$ and $\kappa$ denote the shape operator of $M$ and
 the plane curvature of the plane curve $C$,
respectively.

Since the parametrization $(x(s),0,z(s))$ of the plane curve $C$ is of unit speed,
there exists
a smooth function $\t=\t(s)$ such that $x' = \cos \t$ and $z' = \sin \t$.
Then, the Gaussian curvature $K$ and the mean curvature $H$ of $M$ are, respectively,
 given by
\begin{equation}
\tag{3.4}
\begin{aligned}
K=k_1k_2&=\frac{\t'(s)\sin \t}{x},\\
2H=k_1+k_2&=\t'(s)+\frac{\sin \t}{x}.
\end{aligned}
\end{equation}
Hence, the gradient $\nabla K$ of the Gaussian curvature $K$ of $M$
is given by
\begin{equation}
\tag{3.5}
\begin{aligned}
\nabla K=K'(s)e_1,
\end{aligned}
\end{equation}
where
\begin{equation}
\tag{3.6}
\begin{aligned}
K'(s)=\frac{1}{x^2}\{x\t''(s)\sin \t +x\t'(s)^2\cos \t-\t'(s)\cos\t\sin\t\}.
\end{aligned}
\end{equation}

We now suppose that the Gauss map $G$ of the surface of revolution $M$
satisfies for a $3\times 3$ matrix $A=(a_{ij})$
\begin{equation}
\tag{3.7}
\begin{aligned}
\square G=AG.
\end{aligned}
\end{equation}
Recall that the Gauss map $G$ is given by
\begin{equation}
\tag{3.8}
\begin{aligned}
G(s,t)=(-\sin\t\cos t, -\sin\t\sin t,\cos\t).
\end{aligned}
\end{equation}

Then, it follows from (2.1), (3.4) and (3.5)  that
\begin{equation}
\tag{3.9}
\begin{aligned}
&\{K'(s)\cos\t-2KH\sin\t\}\cos t\\
&=-a_{11}\sin\t\cos t-a_{12}\sin\t\sin t+a_{13}\cos\t,
\end{aligned}
\end{equation}
\begin{equation}
\tag{3.10}
\begin{aligned}
&\{K'(s)\cos\t-2KH\sin\t\}\sin t\\
&=-a_{21}\sin\t\cos t-a_{22}\sin\t\sin t+a_{23}\cos\t
\end{aligned}
\end{equation}
and
\begin{equation}
\tag{3.11}
\begin{aligned}
K'(s)\sin\t+2KH\cos \t=-a_{31}\sin\t\cos t-a_{32}\sin\t\sin t+a_{33}\cos\t.
\end{aligned}
\end{equation}

First, we suppose that the set $J=\{s\in I|\t'(s)\ne 0\}$ is nonempty.
Then $\t(I)$ contains an interval, hence we get from (3.9)-(3.11) that
$a_{12}=a_{13}=a_{21}=a_{23}=a_{31}=a_{32}=0$ and $a_{11}=a_{22}$.
Thus we obtain  $A=diag(\l,\l,\m)$,
\begin{equation}
\tag{3.12}
\begin{aligned}
K'(s)\cos\t-2KH\sin\t=-\l\sin\t,
\end{aligned}
\end{equation}
 and
\begin{equation}
\tag{3.13}
\begin{aligned}
K'(s)\sin\t+2KH\cos \t=\m\cos\t.
\end{aligned}
\end{equation}
Note that (3.12) and (3.13) are equivalent to the following:
\begin{equation}
\tag{3.14}
\begin{aligned}
K'(s)=a\cos\t\sin\t,
\end{aligned}
\end{equation}
 and
\begin{equation}
\tag{3.15}
\begin{aligned}
2KH=-a\sin^2\t+\m,
\end{aligned}
\end{equation}
where we put $a=\m-\l$.

We prove the following lemma, which plays a crucial role in the proof of our main theorem.
\vskip.1in

 \noindent{\bf Lemma 3.1.} Let $M$ be a surface of revolution given by (3.1) with
 nonempty set $J=\{s\in I|\t'(s)\ne 0\}$.
 Suppose that the Gauss map $G$ of $M$  satisfies $\square G = AG$
 for some $3\times3$ matrix $A$.
 Then $A$ is  of the form $\l I$, where $I$ is an identity matrix.

\noindent {\bf Proof.}  The above discussions show that
 $A$ is a diagonal matrix of the form $A=diag(\l,\l,\m)$ for some constants $\l$ and $\m$.
We put $a=\m-\l$.
 Then, it follows from (3.4), (3.6), (3.14) and (3.15) that
 \begin{equation}
\tag{3.16}
\begin{aligned}
x\t''(s)\sin \t +x\t'(s)^2\cos \t-\t'(s)\cos\t\sin\t=ax^2\cos\t\sin\t
\end{aligned}
\end{equation}
and
\begin{equation}
\tag{3.17}
\begin{aligned}
x\t'(s)^2\sin \t +\t'(s)\sin^2 \t=(-a\sin^2\t+\m)x^2.
\end{aligned}
\end{equation}

By differentiating the both sides of (3.17) with respect to $s$, we get
\begin{equation}
\tag{3.18}
\begin{aligned}
&\t''(s)\sin^2 \t+2x\t'(s)\t''(s)\sin \t +x\t'(s)^3\cos\t+3\t'(s)^2\sin \t\cos \t\\
&+2ax^2\t'(s)\sin \t\cos \t=2x\cos\t(-a\sin^2\t+\m).
\end{aligned}
\end{equation}
If we substitute $\t''(s)$ in (3.16) into (3.18), then we have
\begin{equation}
\tag{3.19}
\begin{aligned}
&-x^2\t'(s)^3\cos \t+4x\t'(s)^2\cos\t\sin \t +\{\cos\t\sin^2 \t +4ax^3\cos\t\sin \t\}\t'(s)\\
&+3ax^2\cos\t\sin^2 \t-2\m x^2\cos \t=0.
\end{aligned}
\end{equation}

Let us substitute $\t'(s)^2$ in (3.17) into (3.19). Then we obtain
\begin{equation}
\tag{3.20}
\begin{aligned}
&5x\t'(s)^2\cos\t\sin \t +\{\cos\t\sin^2 \t +5ax^3\cos\t\sin \t-\m x^3\cot\t\}\t'(s)\\
&+3ax^2\cos\t\sin^2 \t-2\m x^2\cos \t=0.
\end{aligned}
\end{equation}
Once more, we substitute $\t'(s)^2$ in (3.17) into (3.20). Then we get
\begin{equation}
\tag{3.21}
\begin{aligned}
\t'(s)=\frac{\g x^2}{\a x^3+\b},
\end{aligned}
\end{equation}
where we put
\begin{equation}
\tag{3.22}
\begin{aligned}
\a(s)&=-5a\sin^2\t(s) +\m,  \b(s)=4\sin^3\t(s),  \\
\g(s)&=-2a\sin^3\t(s)+3\m\sin\t(s).
\end{aligned}
\end{equation}

Now, we replace $\t'(s)$ in (3.17) with that in (3.21). Then we have
\begin{equation}
\tag{3.23}
\begin{aligned}
a_6x^6+a_3x^3+a_0=0,
\end{aligned}
\end{equation}
where we use the following notations:
\begin{equation}
\tag{3.24}
\begin{aligned}
a_6(\t)&=25a^3\cos^6\t+5a^2(15\l-8\m)\cos^4\t\\
&+a(5a-\m)(4\m-15\l)\cos^2\t+\l(5a-\m)^2,
\end{aligned}
\end{equation}
\begin{equation}
\tag{3.25}
\begin{aligned}
a_3(\t)=26a^2\sin^7\t-19a\m \sin^5\t-4\m^2\sin^3\t
\end{aligned}
\end{equation}
and
\begin{equation}
\tag{3.26}
\begin{aligned}
a_0(\t)=-8a\sin^8\t+4\m\sin^6\t.
\end{aligned}
\end{equation}

Let us differentiate (3.23) with respect to $s$.
Here, we denote by $\dot{a_i}(\t)$ the derivative of $a_i(\t)$ with respect to $\t, i=0,3,6$.
Using $x'=\cos \t$ and $\t'(s)$ given by (3.21), we get
\begin{equation}
\tag{3.27}
\begin{aligned}
b_6x^6+b_3x^3+b_0=0,
\end{aligned}
\end{equation}
where we denote
\begin{equation}
\tag{3.28}
\begin{aligned}
b_6(\t)=6\a\cos\t a_6(\t)+\g \dot{a_6}(\t),
\end{aligned}
\end{equation}
\begin{equation}
\tag{3.29}
\begin{aligned}
b_3(\t)=3\a\cos\t a_3(\t)+6\b\cos\t a_6(\t)+\g \dot{a_3}(\t)
\end{aligned}
\end{equation}
and
\begin{equation}
\tag{3.30}
\begin{aligned}
b_0(\t)=3\b\cos\t a_3(\t)+\g \dot{a_0}(\t).
\end{aligned}
\end{equation}

If we compute $b_i(\t)$ for $i=0,3,6$, then we have
\begin{equation}
\tag{3.31}
\begin{aligned}
b_6(\t)=\{1050a^4\sin^8\t+ \sum^6_{i=0}p_i(\l,\m)\sin^i\t\}\cos\t,
\end{aligned}
\end{equation}
\begin{equation}
\tag{3.32}
\begin{aligned}
b_3(\t)=\{-214a^3\sin^9\t+ \sum^7_{i=0}q_i(\l,\m)\sin^i\t\}\cos\t
\end{aligned}
\end{equation}
and
\begin{equation}
\tag{3.33}
\begin{aligned}
b_0(\t)=\{440a^2\sin^{10}\t+ \sum^8_{i=0}r_i(\l,\m)\sin^i\t\}\cos\t,
\end{aligned}
\end{equation}
where $p_i(\l,\m), q_i(\l,\m)$ and $r_i(\l,\m)$ are respectively some polynomials in
$\l$ and $\m$.

Eliminating $x^6$, it follows from (3.23) and (3.27) that
\begin{equation}
\tag{3.34}
\begin{aligned}
c_3x^3+c_0=0,
\end{aligned}
\end{equation}
where
\begin{equation}
\tag{3.35}
\begin{aligned}
c_3=a_3b_6-b_3a_6, c_0=a_0b_6-b_0a_6.
\end{aligned}
\end{equation}
Due to (3.24)-(3.26) and (3.31)-(3.33), we may compute $c_3$ and $c_0$ as follows:
\begin{equation}
\tag{3.36}
\begin{aligned}
c_3=\{32650a^6\sin^{15}\t +\sum^{13}_{j=0}p_{3j}(\l,\m)\sin^j\t\}\cos\t ,
\end{aligned}
\end{equation}
and
\begin{equation}
\tag{3.37}
\begin{aligned}
c_0=\{-19400a^5\sin^{16}\t +\sum^{14}_{j=0}p_{0j}(\l,\m)\sin^j\t\}\cos\t,
\end{aligned}
\end{equation}
where each $p_{ij}(\l,\m)(i=0,3)$ is a polynomial in
$\l$ and $\m$.

Let us replace $x^3$ in (3.23) with $x^3=-c_0/c_3$ given in (3.34).
Then we have
\begin{equation}
\tag{3.38}
\begin{aligned}
a_6c_0^2-a_3c_0c_3+a_0c_3^2=0.
\end{aligned}
\end{equation}
Using (3.31)-(3.33), (3.36) and (3.37), we may compute the leading terms of
those in (3.38)  as follows:
\begin{equation}
\tag{3.39}
\begin{aligned}
a_6c_0^2&=-25(19400)^2a^{13}\sin^{40}\t+ \text{lower degree terms in}  \sin\t,\\
a_3c_0c_3&=26(19400)(32650)a^{13}\sin^{40}\t+\text{lower degree terms in}  \sin\t,\\
a_0c_3^2&=8(32650)^2a^{13}\sin^{40}\t+\text{lower degree terms in}  \sin\t.\\
\end{aligned}
\end{equation}
Hence we obtain
\begin{equation}
\tag{3.40}
\begin{aligned}
a_6c_0^2-a_3c_0c_3+a_0c_3^2&= -17349480000a^{13} \sin^{40}\t\\
&+\text{lower degree terms in}  \sin\t.
\end{aligned}
\end{equation}

Since $\t(I)$ contains an interval,
together with (3.38), (3.40)  shows that $a$ must be zero.
Thus we have $\m=\l$ and hence $A=\l I$. This completes the proof. $\square$

    \vskip.5cm

    \section{Main Theorems and Corollaries}
   \vskip.5cm

Finally, we prove the main theorem as follows.

 \noindent{\bf Theorem 4.1.}
  Let $M$ be a surface of revolution.
 Then  the Gauss map $G$ of $M$  satisfies $\square G = AG$
 for some $3\times3$ matrix $A$ if and only if
  $M$ is an open part of the following surfaces:

\noindent 1)   a plane,

\noindent 2) a right circular cone,

\noindent 3)   a circular cylinder,

\noindent 4) a sphere.

\noindent {\bf Proof.} We consider a surface of revolution $M$ obtained
by  rotating the unit speed plane curve $C:(x(s), 0, z(s))$ with $x(s)>0$ around $z$-axis
which is defined on an interval $I$.

Suppose that the Gauss map $G$ of $M$  satisfies $\square G = AG$
 for some $3\times3$ matrix $A$.
 For a function $\t=\t(s)$ satisfying  $(x'(s),z'(s))=(\cos \t(s),\sin\t(s))$,
 let us put $J=\{s\in I|\t'(s)\ne 0\}$.

 We divide by two cases.

 \noindent {\bf Case 1.} Suppose that $J$ is nonempty. Then, as in the proof of Lemma 3.1
 we have  $A=diag(\l,\l,\m)$  with $a=\l-\m$.
 Furthermore, Lemma 3.1 shows   that $a=0$ that is, $\l=\m$.
  Hence it follows from (3.14) and (3.15) that
  the Gaussian curvature $K$ is constant and the mean curvature $H$
  satisfies $2KH=\l$.

 If $\l\ne 0$, then both of $K$ and $H$ are nonzero constant.
 Hence it follows from a well-known
theorem (cf. \cite{l}) that M is an open part of  a sphere.
Using (3.17) and (3.21) with $a=0$,  it can be directly shown that
$\t'(s)$ is constant and $x(s)=r \sin \t$ for a positive constant $r$.
This shows that the profile curve $C$ is an open part of a half circle centered on the
rotation axis of $M$. Thus, $M$ is an open portion of a round sphere.

If $\l=0$, then $K$ is constant with $2KH=0$.  Suppose that $K\ne 0$.
Then we have $H=0$. But catenoids are the only minimal nonflat surfaces of revolution,
of which  Gaussian curvature $K$ are nonconstant.
This contradiction shows that $K= 0$. Thus $M$ is a flat  surfaces of revolution.
Therefore, $M$ is an open part of a plane, a right circular cone or a circular cylinder.

 \noindent {\bf Case 2.} Suppose that $J$ is empty.  Then the profile curve $C$ of $M$ is a
 straight line. Thus,  $M$ is an open part of a plane, a right circular cone or a circular cylinder.

  The converse is obvious from (2.1). $\square$

Combining the results of \cite{dpv, ks}, the following characterization theorems
can be obtained.
\vskip.3cm

\noindent{\bf Corollary 4.2.}
Let $M$ be a surface of revolution. Then the following are equivalent.

\noindent 1)  $M$ is an open part of a round sphere.

\noindent 2)  The Gauss map $G$ of $M$  satisfies $\square G = AG$
 for some nonsingular $3\times3$ matrix $A$.

\vskip.3cm
\noindent{\bf Corollary 4.3.}
Let $M$ be a surface of revolution. Then the following are equivalent.

\noindent 1)  $M$ is an open part of a right circular cone.

\noindent 2)  The Gauss map $G$ of $M$  satisfies $\square G = AG$
 for some  $3\times3$ matrix $A$, but not satisfies $\Delta G = AG$ for any
  $3\times3$ matrix $A$.


\vskip.8cm
\begin{thebibliography}{10}
\bibitem{ag}
Luis J. Alias and N. Gurbuz, {\itshape An extension of Takahashi theorem for the linearized operators of the higher order mean curvatures}, Geom. Dedicata 121 (2006), 113-127.

\bibitem{bb}
C. Baikoussis and D. E. Blair, {\itshape On the Gauss map of ruled surfaces}, Glasgow Math. J.
34 (1992), no. 3, 355-359.

\bibitem{bv} C. Baikoussis and L. Verstraelen, {\itshape On the Gauss map of helicoidal surfaces}, Rend. Sem.
Mat. Messina Ser. II 2(16) (1993), 31-42.

\bibitem{c1} B.-Y. Chen, {\itshape Total mean curvature and
submanifolds of finite type}, World Scientific
Publ.,  New Jersey (1984).

\bibitem{c2} B.-Y. Chen, {\itshape Finite type
submanifolds and generalizations}, University
of Rome
(1985).
 
 \bibitem{cp} B.-Y. Chen and P. Piccinni, {\itshape  Submanifolds with finite type Gauss map},
  Bull. Austral. Math. Soc. 35 (1987), no. 2, 161-186.

\bibitem{cy}
S. Y. Cheng and  S. T. Yau, {\itshape Hypersurfaces with constant scalar curvature},
 Math. Ann. 225 (1977), no. 3, 195-204.

 \bibitem{cho1}
 S. M. Choi, {\itshape On the Gauss map of surfaces of revolution in a 3-dimensional Minkowski
space}, Tsukuba J. Math. 19 (1995), no. 2, 351-367.

 \bibitem{cho2}  S. M. Choi, {\itshape On the Gauss map of ruled surfaces in a 3-dimensional Minkowski space},
Tsukuba J. Math. 19 (1995), no. 2, 285-304.

\bibitem{ckky}  S. M. Choi, D.-S. Kim, Y. H. Kim and D. W. Yoon, {\itshape Circular cone and its Gauss map},
 Colloq. Math. 129 (2012), no. 2, 203-210.

\bibitem{dpv}
F. Dillen, J. Pas and L.  Verstraelen,  {\itshape On the Gauss map of surfaces of revolution},
 Bull. Inst. Math. Acad. Sinica 18 (1990), no. 3, 239-246.

\bibitem{do}
Manfredo P. do Carmo,  {\itshape Differential geometry of curves and surfaces}. Translated from the Portuguese,
 Prentice-Hall, Inc., Englewood Cliffs, N.J., 1976.

\bibitem{du}
U. Dursun, {\itshape Flat surfaces in the Euclidean space $E^3$   with pointwise 1-type Gauss map},
 Bull. Malays. Math. Sci. Soc. (2) 33 (2010), no. 3, 469-478.
 
\bibitem{kkkr}
 U-H. Ki,  D.-S. Kim, Y. H. Kim and Y.-M. Roh, {\itshape Surfaces of revolution with pointwise 1-type Gauss map in Minkowski 3-space}, Taiwanese J. Math. 13 (2009), no. 1, 317-338.

\bibitem{k1}
 D.-S. Kim, {\itshape On the Gauss map of quadric hypersurfaces},
  J. Korean Math. Soc. 31 (1994), no. 3, 429-437.

  \bibitem{k2}
 D.-S. Kim, {\itshape On the Gauss map of hypersurfaces in the space form},
  J. Korean Math. Soc. 32 (1995), no. 3, 509-518.

\bibitem{kk} D.-S. Kim and Y. H. Kim,
{\itshape Surfaces with planar lines of curvature},  Honam Math. J. 32 (2010), 777-790.

\bibitem{kky}
D.-S. Kim, Y. H. Kim and D. W. Yoon, {\itshape Extended B-scrolls and their Gauss maps},
Indian J. Pure Appl. Math. 33 (2002), no. 7, 1031-1040.

\bibitem{ks}
 D.-S. Kim and B. Song, {\itshape On the Gauss map of generalized slant
       cylindrical surfaces}, J. Korea Soc. Math. Educ. Ser. B: Pure Appl. Math.,
       20(3), 149-158.

\bibitem{kt}
Y. H. Kim and N. C. Turgay, {\itshape Surfaces in $E^3$   with $L_1$-pointwise 1-type Gauss map},
 Bull. Korean Math. Soc. 50 (2013), no. 3, 935-949.

\bibitem{ky}
Y. H. Kim and D. W.  Yoon, {\itshape On the Gauss map of ruled surfaces in Minkowski space},
  Rocky Mountain J. Math. 35 (2005), no. 5, 1555-1581.

\bibitem{l}
T. Levi-Civita, {\itshape Famiglie di superficie isoparametriche nell¡¯ordinario spacio euclideo},
 Rend. Acad. Lincei 26 (1937), 355-362.

\bibitem{rv}
E. A. Ruh and J. Vilms, {\itshape The tension field of the Gauss map}, Trans. Amer. Math. Soc.
149 (1970), 569-573.
\end{thebibliography}
\end{document}